\documentclass[a4paper,10pt]{scrartcl}
\usepackage{amsmath,amsthm}
\usepackage{amsfonts}
\usepackage{amssymb}
\usepackage{a4wide}
\usepackage{amscd}
\usepackage{color}
\usepackage{graphicx}
\usepackage{float}
\usepackage[utf8]{inputenc}
\usepackage{cleveref}
\usepackage{algorithm}
\usepackage{algpseudocode}

\usepackage{multirow}

\newcommand{\N}{\mathbb{N}}
\newcommand{\R}{\mathbb{R}}

\newcommand{\norm}[2]{\| #1 \|_{#2}}

\newcommand{\yd}{y^\delta}
\newcommand{\xdag}{x^\dag}
\newcommand{\xad}{x_{\alpha}^\delta}

\newcommand{\supp}{\operatorname{supp}}

\DeclareMathOperator*{\argmin}{argmin}

\newenvironment{frcseries}{\fontfamily{frc}\selectfont}{}

\newtheoremstyle{slanted}
{}
{}
{\slshape}
{}
{\bfseries}
{.}
{ }
{}
\newtheoremstyle{roman}{}{}{\rmfamily}{}{\bfseries}{.}{ }{}

\theoremstyle{plain}
\newtheorem{theorem}{Theorem}
\newtheorem{definition}{Definition}
\newtheorem{lemma}{Lemma}
\newtheorem{proposition}{Proposition}
\theoremstyle{slanted}

\theoremstyle{roman}

\title{Deautoconvolution in the two-dimensional case}
\author{\textsc{Yu Deng}\thanks{Chemnitz University of Technology, Faculty of Mathematics, 09107 Chemnitz, Germany, e-mail: yu.deng@math.tu-chemnitz.de},\,\,
\textsc{Bernd Hofmann}\thanks{Faculty of Mathematics, Chemnitz University of Technology, 09107 Chemnitz, Germany, \newline e-mail: hofmannb@mathematik.tu-chemnitz.de.}\,\, and \textsc{Frank Werner}\thanks{Institute for Mathematics, University of W\"urzburg, Emil-Fischer-Str.~30, 97074~W\"urzburg, Germany, e-mail: frank.werner@mathematik.uni-wuerzburg.de}}
\date{\today}

\begin{document}

\maketitle

\begin{quote}
\textbf{Abstract:} There is extensive mathematical literature on the inverse problem of deautoconvolution for a function with support in the unit interval $[0,1] \subset \R$, but little is known about the multidimensional situation. This article tries to fill this gap with analytical and numerical studies on the reconstruction of a real function of two real variables over the unit square from observations of its autoconvolution on $[0,2]^2 \subset \R^2$ (full data case) or on $[0,1]^2$ (limited data case). In an $L^2$-setting, twofoldness and uniqueness assertions are proven for the deautoconvolution problem in 2D. Moreover, its ill-posedness is characterized and illustrated. Extensive numerical
case studies give an overview of the behaviour of stable approximate solutions to the two-dimensional deautoconvolution problem obtained by Tikhonov-type regularization  with different penalties and the iteratively regularized Gauss-Newton method.
\end{quote}

\begin{quote}
\end{quote}

\begin{quote}
\noindent
{\small \textbf{AMS-classification (2010):}
47J06, 65R32, 45Q05, 47A52, 65J20
}
\end{quote}
\section{Introduction} \label{sec:intro}

The object of research in this work is the problem of \textsl{deautoconvolution}, where our focus is on the \textsl{two-dimensional case}, which means that a square integrable real function of two variables $x(t_1,t_2)\;(0 \le t_1,t_2 \le 1)$
is to be identified from the function $y=x*x$ of its autoconvolution. If we consider $x$ as an element in the Hilbert space $L^2(\R^2)$ with support $\supp(x) \subseteq [0,1]^2$, then it is well-known that $x*x$ also lies in $L^2(\R^2)$ with support $\supp(x*x) \subseteq [0,2]^2$. In this context, the elements $x$ and $x*x$ both can be considered as tempered distributions with compact support, where $\supp(\cdot)$ is regarded as the essential support with respect to the Lebesgue measure $\lambda$ in $\R^2$. Instead of $y$ itself, noisy data $\yd \in  L^2(\R^2)$ to $y$ with some noise level $\delta \ge 0$ are available only. Since the inverse problem of deautoconvolution tends to be ill-posed, the aim of the recovery process
is to find \textsl{stable approximate solutions} of $x$ based on the data $\yd$.
We are going to distinguish the \textit{full data case}, where noisy data are available for $y(s_1,s_2)\; (0 \le s_1,s_2 \le 2)$, and the \textsl{limited data case}, where data are given for $y(s_1,s_2)\; (0 \le s_1,s_2 \le 1)$.
Since the scope of the data in the limited data case is only 25\% compared to the full data case, the effect of ill-posedness is stronger in that case.
As a consequence, also the chances for the accurate recovery of $x$ are more restricted in the limited data case.

The simplest application of our deautoconvolution problem in two dimensions is the recovery of the density function $x$ with support in the unit square $[0,1]^2$ of a two-dimensional random variable $\mathfrak{X}$ from observations
of the density function $y$ of the two-dimensional random variable $\mathfrak{Y}:=\hat{\mathfrak{X}}+\bar{\mathfrak{X}}$, where $\mathfrak{X}, \hat{\mathfrak{X}}$ and $\bar{\mathfrak{X}}$ are assumed to be of i.i.d.~type.

The deautoconvolution problem \textsl{in one dimension} has been considered extensively in the literature motivated by physical applications in spectroscopy (see, e.g., \cite{Baumeister91,Schleicher83}). Its mathematical analysis has been implemented comprehensively in the last decades with focus on properties of the specific forward operator, ill-posedenss and regularization based on the seminal paper \cite{GorHof94}. In this context, we refer to
\cite{BuerFlem15,BuerHof15,BuerMat17,ChoiLan05,FleiGorHof99,FleiHof96,Janno00} for investigations concerning the stable identification of real functions $x$ on the unit interval $[0,1]$ from noisy data of its autoconvolution $x*x$. A new series of
interdisciplinary autoconvolution studies was unleashed by a cooperation started in 2010 between a Research Group of the Max Born Institute for Nonlinear Optics and Short Pulse Spectroscopy, Berlin, led by Prof.~G\"unter Steinmeyer and the Chemnitz research Group on Regularization,
and we refer to the publications \cite{Anzengruber16,BuerFlemHof16,Flemming18,Gerth14} presenting the output of this cooperation. The goal of this cooperation between mathematics and laser optics was the extension of the one-dimensional deautoconvolution problem to complex-valued functions combining amplitude and phase functions for characterizing ultrashort laser pulses.

In this article, in an $L^2$-setting we are considering
a series of numerical case studies for the \textsl{nonlinear Volterra-type integral equation}
\begin{equation}\label{eq:auto}
	F(x)=y\,,\quad \mbox{with} \quad F(x):=x*x\,,
\end{equation}
the solution of which solves the \textsl{deautoconvolution problem in two dimensions}. Equation \eqref{eq:auto} is a special case of a nonlinear operator equation
\begin{equation}\label{eq:opeq}
	F(x)=y,\quad F:\mathcal{D}(F) \subseteq X \to Y\,,
\end{equation}
with \textsl{forward operator} $F$ mapping between \textsl{real-valued Hilbert spaces} $X$ and $Y$ with norms $\|\cdot\|_X$ and $\|\cdot\|_Y$, respectively, and domain $\mathcal{D}(F)$.

In dependence of the data situation, we have to distinguish  in the full data case the forward operator  $F: X=L^2([0,1]^2) \to  Y=L^2([0,2]^2)$ defined as
\begin{small}
	\begin{equation}\label{eq:auto1}
		[F(x)](s_1,s_2):=\!\! \! \!\int \limits_{\max(s_2-1,0)}^{\min(s_2,1)}\int \limits_{\max(s_1-1,0)}^{\min(s_1,1)} \! \! x(s_1-t_1,s_2-t_2)\,x(t_1,t_2)\,d t_1 \, d t_2 \quad  (0 \le s_1, s_2 \le 2)
	\end{equation}
\end{small}
and in the limited data case the forward operator $F: X=L^2([0,1]^2) \to  Y=L^2([0,1]^2)$ as
\begin{equation}\label{eq:auto2}
	[F(x)](s_1,s_2):=\;\; \int _0^{s_2} \int_0^{s_1} x(s_1-t_1,s_2-t_2)\,x(t_1,t_2)\,d t_1\, d t_2 \quad  (0 \le s_1, s_2 \le 1).
\end{equation}
In general we consider $\mathcal{D}(F)=X=L^2([0,1]^2)$, but for the limited data case we partially focus on \textsl{non-negative solutions} expressed by the domain $\mathcal{D}(F)=\mathcal{D}^+$ with
\begin{equation} \label{eq:domainnonnegative}
	\mathcal{D}^+:=\{x \in X=L^2([0,1]^2):\, x \ge 0 \;\; \mbox{a.e.~on}\;[0,1]^2\}\,.
\end{equation}

For any function $x \in L^2([0,1]^2)$ the autoconvolution products $F(x)=x * x$ and $F(-x)=(-x)*(-x)$ coincide for both forward operator versions \eqref{eq:auto1} and \eqref{eq:auto2}.
However, it is of interest whether for $y=x*x$ the elements $x$ and $-x$ are the only solutions of equation \eqref{eq:auto} or not.
Moreover it is of interest whether in the limited data case the restriction of the domain $\mathcal{D}(F)$ to $\mathcal{D}^+$ from \eqref{eq:domainnonnegative} leads to unique solutions.
Some answers to those questions will be given in the subsequent Section~\ref{sec:unique}.

\medskip

The remainder of the paper is organized as follows: Section~\ref{sec:unique} is devoted to assertions on twofoldness and uniqueness for the deautoconvolution problem in two dimensions, preceded by a subsection with relevant
lemmas and definitions. As an inverse problem, deautoconvolution tends to be ill-posed in the setting of infinite dimensional $L^2$-spaces.
After the presentation of two functions defined over the unit square as basis for later numerical case studies, in Section~\ref{sec:ill-posed} the specific ill-posedness character for the deautoconvolution of a real function of two real variables with compact support is analyzed and illustrated. To suppress ill-posedness phenomena, variants of variational and iterative regularization methods are used, which will be introduced in Section~\ref{sec:regu}.
The numerical treatment, including discretization approaches of forward operator and penalty functionals for the Tikhonov regularization as well as for the iterative  regularization by using the Fourier transform, is outlined in Section~\ref{sec:numerics}. Section~\ref{sec:studies} completes the article with comprehensive
numerical case studies.

\section{Assertions on twofoldness and uniqueness for the deautoconvolution problem in two dimensions}\label{sec:unique}

\subsection{Preliminaries}

Assertions on twofoldness and uniqueness for the deautoconvolution problem in one dimension have been formulated in the articles \cite{GorHof94} for the limited data case and \cite{Gerth14} for the full data case. The respective proofs
are based on the Titchmarsh convolution theorem from \cite{Titchmarsh26}, which was formulated as Lemma~3 in \cite{GorHof94} and will be recalled below in a slightly reformulated form as Lemma~\ref{lem:tit1}.

\begin{lemma} \label{lem:tit1}
	Let the functions $f,g \in L^2(\R)$ have compact supports $\supp(f)$ and $\supp(g)$. Then we have for the convolution that $f*g \in L^2(\R)$ and that the equation
	\begin{equation} \label{eq:sum1}
		\inf \supp (f*g)\,=\, \inf \supp(f)+ \inf \supp(g)
	\end{equation}
	holds. In particular, for $\supp(f)$ and $\supp(g)$ covered by the unit interval $[0,1]$, we conclude from
	$$[f*g](s)=\int \limits_{\max(s-1,0)}^{\min(s,1)} f(s-t)\,g(t)\,dt\,=\, 0 \quad \mbox{a.e.\, for}\;\; s \in [0,\gamma] \quad(\gamma \le 2)$$
	that there are numbers $\gamma_1,\gamma_2 \in [0,1]$ with $\gamma_1+\gamma_2 \ge \gamma$ such that
	$$ f(t)\, = \, 0 \quad \mbox{a.e.\, for}\;\; t \in [0,\gamma_1] \quad \mbox{and} \quad g(t) \,=\, 0 \quad \mbox{a.e.\, for}\;\; t \in [0,\gamma_2].  $$
\end{lemma}

For an extension of the Titchmarsh convolution theorem  to two dimensions, we mention the following Lemma~\ref{lem:tit2} (cf.~\cite{Lions51,Lions53}).

\begin{lemma} \label{lem:tit2}
	Let the functions $f,g \in L^2(\R^2)$ have compact supports $\supp(f)$ and $\supp(g)$. Then we have for the convolution that $f*g \in L^2(\R^2)$ and that the equation
	\begin{equation} \label{eq:sum2}
		\mathrm{ch\,supp}(f*g)\,=\, \mathrm{ch\;supp}(f)+\mathrm{ch\,supp}(g)
	\end{equation}
	holds, where $\mathrm{ch}\; M$ denotes the convex hull of a set $M \subseteq \R^2$. In the special case that $\supp(f*g)=\emptyset$ we have that at least one of the supports $\supp(f)$ or $\supp(g)$ is the empty set.
\end{lemma}

\begin{definition} \label{def:full}
	For given $y \in L^2([0,2]^2)$, we call $\xdag \in L^2([0,1]^2)$ with $\supp(\xdag) \subseteq [0,1]^2$ a \textsl{solution} to the operator equation \eqref{eq:auto} \textsl{in the full data case}
	if it satisfies the condition
	\begin{equation} \label{eq:full}
		[\xdag * \xdag](s_1,s_2)=y(s_1,s_2) \quad \mbox{a.e. for} \quad  (s_1,s_2) \in [0,2]^2\,.
	\end{equation}
\end{definition}

\begin{definition} \label{def:limited}
	For given $y \in L^2([0,1]^2)$, we call $\xdag \in L^2([0,1]^2)$ with $\supp(\xdag) \subseteq [0,1]^2$ a \textsl{solution} to the operator equation \eqref{eq:auto} \textsl{in the limited data case}
	if it satisfies the condition
	\begin{equation} \label{eq:limited}
		[\xdag * \xdag](s_1,s_2)=y(s_1,s_2) \quad \mbox{a.e. for} \quad  (s_1,s_2) \in [0,1]^2.
	\end{equation}
	For $\xdag \in \mathcal{D}^+$ with $\mathcal{D}^+$ from \eqref{eq:domainnonnegative} we call it \textsl{non-negative solution in the limited data case}.
	
\end{definition}

\begin{definition} \label{def:factored}
	We call $x \in L^2([0,1]^2)$ with $\supp(x) \subseteq [0,1]^2$ satisfying \eqref{eq:full} or \eqref{eq:limited} a \textrm{factored solution} to equation~(\ref{eq:auto}) in the full data case or in the limited data case, respectively, if we have the structure
	$x(t_1,t_2)=x_1(t_1)\,x_2(t_2)\;(0 \le t_1,t_2 \le 1)$ with $x_i \in L^2([0,1])$, $\supp(x_i)\subseteq [0,1]$ for $i=1$ and $i=2$. If moreover $x_i \ge 0$ a.e. on $[0,1]$ for $i=1$ and $i=2$, then we call it \textrm{non-negative factored solution} in the respective case.
\end{definition}

\subsection{Results for the full data case}

Lemma~\ref{lem:tit2} allows us to prove the following theorem for the forward autoconvolution operator $F: L^2([0,1]^2) \to  L^2([0,2]^2)$ from \eqref{eq:auto1}, which is an extension of \cite[Theorem~4.2]{Gerth14} to the two-dimensional case of the deautoconvolution problem.

\begin{theorem} \label{thm:tit1}
	If, for given $y \in L^2([0,2]^2)$, the function $\xdag \in L^2([0,1]^2)$ with $\supp(\xdag) \subseteq [0,1]^2$  is a solution to (\ref{eq:auto}) with $F$ from \eqref{eq:auto1},
	then $\xdag$ and $-\xdag$ are the only solutions of this equation in the full data case.
\end{theorem}
\begin{proof}
	Let $\xdag \in L^2([0,1]^2)$ supposing $\supp(\xdag) \subseteq [0,1]^2$ and $h \in L^2([0,1]^2)$ supposing $\supp(h) \subseteq [0,1]^2$. We assume that $\xdag$ and $\xdag+h$  solve the equation~(\ref{eq:auto}), which means that $[\xdag*\xdag](s_1,s_2)=[(\xdag+h)*(\xdag+h)](s_1,s_2)$ a.e. for $(s_1,s_2) \in [0,2]^2$. Then we have $[(\xdag+h)*(\xdag+h)-\xdag*\xdag](s_1,s_2)=[h*(2\xdag+h)](s_1,s_2)=0$ a.e. for $(s_1,s_2) \in [0,2]^2$. By setting $f:=h$ and $g:=2\xdag+h$ we now apply Lemma~\ref{lem:tit2}. Taking into account that \linebreak $\supp(h*(2\xdag+h)) \subseteq[0,2]^2$ we then have $\supp(h*(2\xdag+h))=\emptyset$ and consequently also \linebreak $\mathrm{ch\,supp}(h*(2\xdag+h))=\emptyset$. This implies due to equation \eqref{eq:sum2} that either $\supp(h)=\emptyset$ or $\supp(2\xdag+h)=\emptyset$.
	On the one hand, $\supp(h)=\emptyset$ leads to the solution $\xdag$ itself, whereas on the other hand $\supp(2\xdag+h)=\emptyset$ leads to $[2\xdag+h](t_1,t_2)=0$ a.e. for $(t_1,t_2) \in [0,1]^2$
	and consequently with $h=-2\xdag$ to the second solution $-\xdag$. Alternative solutions are thus excluded, which proves the theorem.
\end{proof}

\subsection{Results for the limited data case}

For solutions  $\xdag \in L^2([0,1]^2)$ to equation~(\ref{eq:auto}) with $\supp(\xdag) \subseteq [0,1]^2$,  the condition $0 \in \supp(\xdag)$ plays a prominent role in the limited data case. This condition means that
for any ball $B_r(0)$ around the origin with arbitrary small radius $r>0$ there exists a set $M_r \subset B_r(0) \cap [0,1]^2$ with Lebesgue measure $\lambda(M_r)>0$ such that $\xdag(t_1,t_2) \not=0$ a.e. for $(t_1,t_2) \in M_r$.
Vice versa, for $0 \notin \supp(\xdag)$ we have some sufficiently small radius $r>0$ such that $\xdag(t_1,t_2)=0$ a.e. for $(t_1,t_2) \in B_r(0) \cap [0,1]^2$.

First, we generalize in Theorem~\ref{thm:nonunique} those aspects of \cite[Theorem~1]{GorHof94} that concern the strong non-injectivity of the autoconvolution operator in the limited data case.

\begin{theorem} \label{thm:nonunique}
	If, for given $y \in L^2([0,1]^2)$, the function $\xdag \in L^2([0,1]^2)$ with $\supp(\xdag) \subseteq [0,1]^2$ is a solution to (\ref{eq:auto}) with $F$ from \eqref{eq:auto2} that fulfils the condition
	\begin{equation} \label{eq:notzero1}
		0 \notin \supp(\xdag)\,,
	\end{equation}
	then there exist infinitely many other solutions ${\hat x}^\dagger \in L^2([0,1]^2)$ to (\ref{eq:auto}) with $\supp({\hat x}^\dagger) \subseteq [0,1]^2$ in the limited data case.
\end{theorem}
\begin{proof}
	If \eqref{eq:notzero1} holds, there is some $0<\varepsilon<1/2$ such that $\xdag(t_1,t_2)=0$ a.e. for $(t_1,t_2) \in [0,\varepsilon]^2$.
	Then we have, for all elements $h \in L^2([0,1]^2)$ with $\supp(h) \subseteq [0,1]^2$ satisfying the condition
	$$h(t_1,t_2)\,=\,0 \quad \mbox{a.e. for} \;\; (t_1,t_2) \in [0,1]^2 \setminus [1-\varepsilon,1]^2\,,$$
	that ${\hat x}^\dagger=\xdag+h$ obeys the condition
	$$[{\hat x}^\dagger * {\hat x}^\dagger](s_1,s_2)=y(s_1,s_2) \quad \mbox{a.e. for} \quad  (s_1,s_2) \in [0,1]^2\,.$$
	This is a consequence of the fact that $[h * (2\xdag+h)](s_1,s_2)=0$ a.e. for $(s_1,s_2) \in [0,1]^2$ is true for each such element $h$.
\end{proof}

To formulate uniqueness assertions for solutions $\xdag$  to equation~(\ref{eq:auto}) in the limited data case, we restrict our considerations now to non-negative solutions and the domain $\mathcal{D}(F)=\mathcal{D}^+$
from \eqref{eq:domainnonnegative} for the forward operator $F$ from \eqref{eq:auto2}. We present in Theorem~\ref{thm:factored} a result that extends to the two-dimensional autoconvolution operator $F:\mathcal{D}^+ \subset L^2([0,1]^2) \to  L^2([0,1]^2)$ from \eqref{eq:auto2} those aspects of \cite[Theorem~1]{GorHof94} which concern the solution uniqueness. Precisely, we are able to handle the special case of \textsl{factored non-negative solutions} in the sense of Definition~\ref{def:factored}, occurring for example when $\xdag$ is a density function for the two-dimensional random variable $\mathfrak{X}=(\mathfrak{X}_1,\mathfrak{X}_2)$, where $\mathfrak{X}_1$ and $\mathfrak{X}_2$ are uncorrelated one-dimensional random variables.

\begin{theorem} \label{thm:factored}
	Let, for given $y \in L^2([0,1]^2)$, $\xdag$ be a non-negative factored solution to equation~(\ref{eq:auto}) in the limited data case, which satisfies the condition
	\begin{equation} \label{eq:zero3}
		0 \in \supp(\xdag)\,.
	\end{equation}
	Then there are no other non-negative factored solutions in this case.
\end{theorem}
\begin{proof}
	For the factored situation, we have that the right-hand side $y$ is also factored with $$y(s_1,s_2)=y_1(s_1)\, y_2(s_2)\;(0 \le s_1,s_2 \le 1) \quad \mbox{and} \quad  y_1=\xdag_1*\xdag_1,\;\;y_2=\xdag_2*\xdag_2.$$
	Moreover, the condition \eqref{eq:zero3} implies that
	\begin{equation} \label{eq:implicond}
		\inf \supp(\xdag_1)=\inf \supp(\xdag_2)=0\,.
	\end{equation}
	Otherwise, there would be a square $[0,\varepsilon]^2$ with $\varepsilon=\max\{\inf \supp(\xdag_1),\inf \supp(\xdag_2)\}>0$ on which $\xdag$ vanishes almost everywhere such that $0 \notin \supp(\xdag)$ would then apply.
	Now we suppose that for $i=1$ and $i=2$ quadratically integrable perturbations $h_i(t_i)\;(0 \le t_i \le 1)$ exist such that $\xdag_i+h_i \ge 0$ a.e. on $[0,1]$ and
	\begin{equation} \label{eq:added}
		(\xdag_1+h_1)*(\xdag_1+h_1)=y_1 \quad \mbox{and} \quad (\xdag_2+h_2)*(\xdag_2+h_2)=y_2.
	\end{equation}
	To complete the proof of the theorem we still show that $h_1$ and $h_2$ have to vanish almost everywhere on $[0,1]$. This can be done with the help of Titchmarsh's convolution theorem in the one-dimensional case (cf.~Lemma~\ref{lem:tit1}). From \eqref{eq:added} we derive for $i=1$ and $i=2$ that
	$$[h_i*(2\xdag_i+h_i)](s_i)=0\quad \mbox{a.e. for} \;\; s_i \in [0,1],$$
	where $\xdag_i+h_i \ge 0$ implies that $2\xdag_i+h_i \ge \xdag_i$ and $\inf \supp(2\xdag_i+h_i)=0$ as a consequence of \eqref{eq:implicond}. Then it follows from Lemma~\ref{lem:tit1} that
	$$\inf \supp(h_i)+\inf \supp(2\xdag_i+h_i)=\inf\supp(h_i*(2\xdag_i+h_i))  \ge 1$$ and hence $\inf \supp(h_i) \ge 1$ for both $i=1,2$, which gives $h_i=0$ a.e.~on $[0,1]$ and completes the proof.
\end{proof}

\section{Examples and ill-posedness phenomena of deautoconvolution in 2D} \label{sec:ill-posed}
\subsection{Two Examples}\label{sec:example}
\begin{figure}[H]
	\begin{center}
		\includegraphics[width=4.2cm]{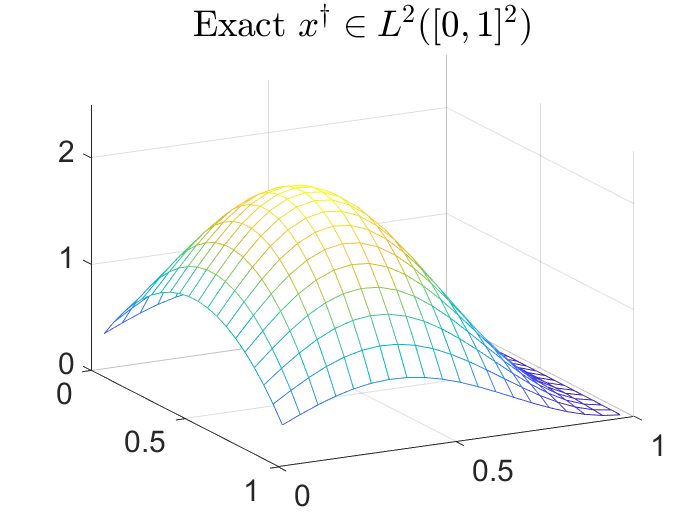}
		\includegraphics[width=4.2cm]{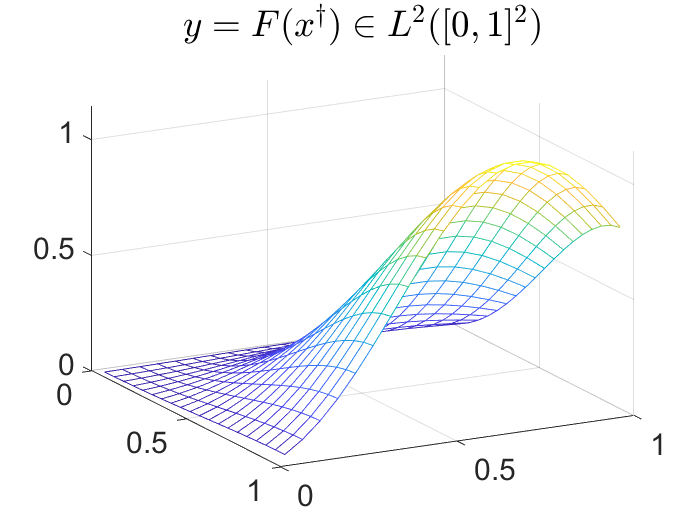}
		\includegraphics[width=4.2cm]{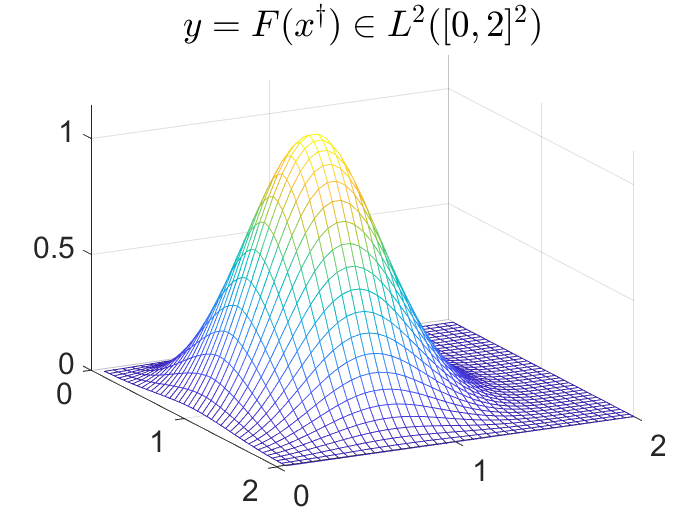}
		\caption{Smooth factored function $x^\dagger(t_1,t_2)$ and $F(x^\dagger)$ with limited data and full data, respectively.}
		\label{fig:contposx}
	\end{center}
\end{figure}

For the numerical case studies of deautoconvolution in 2D, we have selected two examples of solutions $\xdag$ to the autoconvolution equation in 2D. The first example refers to the function
\begin{equation}\label{eq:ex1}
	x^\dagger(t_1,t_2)=\left(-3t_1^2+3t_1+\tfrac14\right)\left(\sin(1.5\pi t_2)+1\right) \quad (0 \le t_1, t_2 \le 1)
\end{equation}
to be reconstructed from its own autoconvolution $F(x^\dagger)=x^\dagger * x^\dagger$. This \textsl{smooth and non-negative factored} function $\xdag$ is illustrated in Figure~\ref{fig:contposx}, in a line with the
$F(\xdag)$-images for the limited and full data case, respectively.

\begin{figure}[H]
	\begin{center}
		\includegraphics[width=4.2cm]{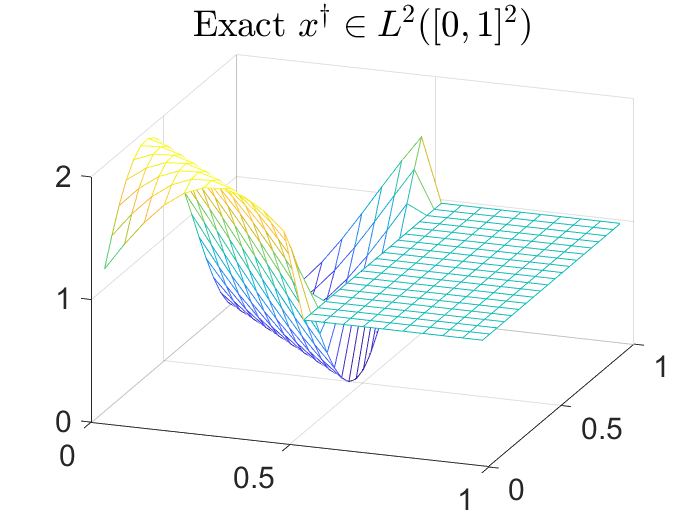}
		\includegraphics[width=4.2cm]{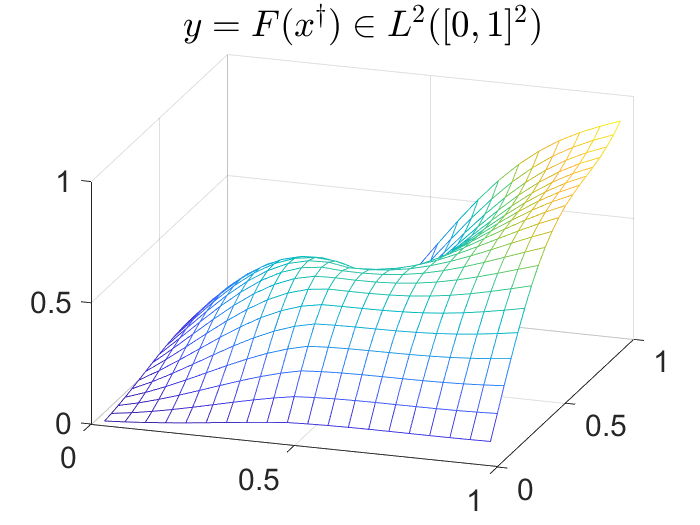}
		\includegraphics[width=4.2cm]{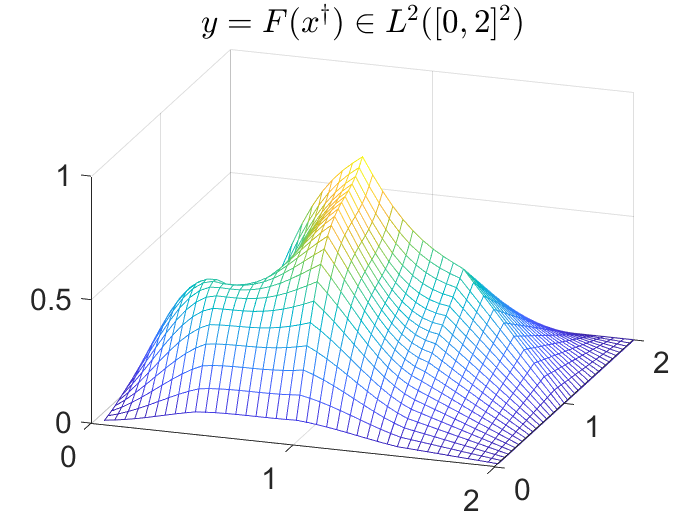}
		\caption{Non-smooth non-factored function $x^\dagger(t_1,t_2)$ and $F(x^\dagger)$ with limited data and full data, respectively.}
		\label{fig:MischPos}
	\end{center}
\end{figure}

The second example refers to the \textsl{non-smooth,  non-factored} and non-negative solution
\begin{equation}\label{eq:ex2}
	x^\dagger(t_1,t_2)=\left\{\begin{array}{ll}
		\sin(1.5\pi (t_1+t_2))+1&\quad (0\leq t_1\leq 0.5, \;0 \le t_2 \le 1)\\
		1&\quad (0.5< t_1\leq 1,  \;0 \le t_2 \le 1)\\
	\end{array}
	\right.,
\end{equation}
which is illustrated in \Cref{fig:MischPos}, again in a line with $F(\xdag)$-images for the limited and full data case, respectively.

\subsection{Ill-posedness phenomenon}
As we will see in particular for the case of limited data in the numerical case studies to be presented below, least-squares solutions to the discretized two-dimensional problem of deautoconvolution tend to become strongly oscillating even if the noise level $\delta>0$ in the observed right-hand side $y^\delta$ is small. This indicates
instability as a phenomenon of \textsl{ill-conditioning} for a discretized variant of deautoconvolution and of \textsl{ill-posedness} for the underlying operator equation \eqref{eq:auto} in infinite dimensional $L^2$-spaces.
For a theoretical verification we adopt the \textsl{concept of local ill-posedness} along the lines of \cite[Def.~1.1]{HofSch98} for nonlinear operator equations, and we recall this concept by the following definition.

\begin{definition}  \label{def:posed}
	An operator equation $F(x)=y$ with nonlinear forward operator $F:\mathcal{D}(F) \subseteq X \to Y$  mapping between the Hilbert spaces $X$ and $Y$ with domain $\mathcal D(F)$ is called \textsl{locally ill-posed} at a solution point
	$\xdag \in \mathcal{D}(F)$ if there exist for
	all closed balls $\overline{\mathcal{B}_r(\xdag)}$ with radius $r>0$ and center $\xdag$
	sequences $\{x_n\} \subset \overline{\mathcal{B}_r(x^\dagger)} \cap \mathcal{D}(F)$ satisfying the condition
	$$\|F(x_n)-F(x^\dagger)\|_Y \to 0 \, , \quad \mbox{but} \quad
	\|x_n-x^\dagger\|_X \not \to 0 \, ,
	\quad \mbox{as} \quad n \to \infty\;.$$
	Otherwise, the operator equation is called  \textsl{locally well-posed} at $\xdag$.
\end{definition}

\smallskip

\textsl{Local ill-posedness everywhere} on the non-negativity domain $$\mathcal{D}(F)=\{x \in X=L^2([0,1]):\, x \ge 0 \;\; \mbox{a.e.~on}\;[0,1]\}$$    was proven for the \textsl{one-dimensional} deautoconvolution problem in the
limited data case in \cite[Lemma~6]{GorHof94}. With the following proposition we extend by using similar proof ideas this assertion to the two-dimensional case and $\mathcal{D}(F)=\mathcal{D}^+$  with $\mathcal{D}^+$ from \eqref{eq:domainnonnegative}.
We should mention as an overall consequence of the observed ill-posedness that the \textit{stable approximate solution} of the two-dimensional decautoconvolution problem \textsl{requires} the use of variational or iterative \textsl{regularization methods}.

\begin{proposition}\label{pro:ill_limited}
	For the limited data case, the operator equation \eqref{eq:auto} with $X=Y=L^2([0,1]^2)$ and forward operator $F$ from \eqref{eq:auto2} restricted to the non-negativity domain $\mathcal{D}(F)=\mathcal{D}^+$ from \eqref{eq:domainnonnegative} is
	\textsl{locally ill-posed everywhere} on $\mathcal{D}^+$.
\end{proposition}
\begin{proof}
	Let $\xdag \in \mathcal{D}^+$ be a solution to the operator equation under consideration here. To show local ill-posedness at $\xdag$ we introduce for fixed $r>0$ a sequence $\{h_n\}_{n=3}^\infty$ of perturbations of the form
	$$h_n(t_1,t_2):=\left\{\begin{array}{ccl} n \, r & \mbox{for} & (t_1,t_2) \in [1-\frac{1}{n},1]^2\\ & & \\0 & \mbox{for} & (t_1,t_2) \in [0,1]^2 \setminus [1-\frac{1}{n},1]^2   \end{array} \right. $$
	with $x_n:=\xdag+h_n \in \mathcal{D}^+$, $\|h_n\|_{L^2([0,1]^2)}=r$ and consequently $x_n \in \overline{{\mathcal B}_r(x^\dagger)} \cap \mathcal{D}^+$ for all $n \ge 3$. To complete the proof of the proposition we still need to show that
	the norm \linebreak $\|F(x_n)-F(\xdag)\|_{L^2([0,1]^2)}$ tends to zero as $n$ tends to infinity. Due to the facts that \linebreak $F(x_n)-F(\xdag)=2\xdag*h_n+h_n*h_n$ and $\|h_n*h_n\|_{L^2([0,1]^2)}=0\,$, it is sufficient to show the limit condition
	$$\|\xdag*h_n\|_{L^2([0,1]^2)}  \to 0 \quad \mbox{as} \quad n \to \infty\,. $$
	Evidently, the non-negative values
	$$ [\xdag*h_n](s_1,s_2)=\int_0^{s_2}\int_0^{s_1} h_n(s_1-t_1,s_2-t_2)\,\xdag(t_1,t_2)\, d t_1 \, d t_2   $$
	can be different from zero only for the pairs $(s_1,s_2) \in [1-\frac{1}{n},1]^2$. Using the Cauchy-Schwarz inequality and taking into account that $\xdag \in \mathcal{D}^+$
	we have for those pairs the estimate
	$$ [\xdag*h_n](s_1,s_2)=  n \, r \int \limits_0^{s_2-(1-\frac{1}{n})}\int \limits_0^{s_1-(1-\frac{1}{n})} \xdag(t_1,t_2)\, d t_1 \, d t_2 \le r\,\|\xdag\|_{L^2([0,1]^2)}.   $$
	This, however, yields $$\|\xdag*h_n\|_{L^2([0,1]^2)} \le r\,\,\|\xdag\|_{L^2([0,1]^2)} \left(\int \limits_{1-\frac{1}{n}}^1 \int \limits_{1-\frac{1}{n}}^1 \, d s_1 \, d s_2\right)^{1/2}=\frac{ r\,\,\|\xdag\|_{L^2([0,1]^2)}}{n}$$
	tending to zero as $n$ tends to infinity, which completes the proof.
\end{proof}

In the full data case of one-dimensional deautocovolution, local ill-posedness everywhere has been shown by Proposition~2.3 in \cite{FleiHof96}. The used counterexample, however, is much more sophisticated and requires perturbations
with weak poles at the origin. This seems to indicate the significantly lower strength of ill-posedness for the full data case compared to the limited data case. For factored solutions, the counterexample from \cite{FleiHof96} can also be exploited to prove local ill-posedness for the two-dimensional deautoconvolution problem in the full data case. Numerical case studies confirm the lower level of instability in the full data case of the 2D deautoconvolution
compared to the limited data case, and we refer to Figure~\ref{fig:illpose} below.

Taking into account throughout the paper the \textsl{deterministic noise model}
\begin{equation}\label{eq:noisylevel}
	\|y^\delta-y\|_Y\leq \delta\,,
\end{equation}
we have calculated \textsl{discretized least-square solutions} $x_{ls}$ of the deautoconvolution problem for  the first example $x^\dagger$ given in \eqref{eq:ex1}. In this context, a discretization with $21 \times 21$ grid points over the
unit square has been exploited for minimizing the Euclidean norm squares of the discretized residuals $F(\xdag)-y^\delta$. By using a noise level $\delta>0$ that corresponds to a relative data error of $0.8\%$, the least squares situations
for $Y:=L^2([0,2]^2)$ of the full data case with $F$ defined in \eqref{eq:auto2} and for $Y:=L^2([0,1]^2)$ of the limited data case with $F$ defined in \eqref{eq:auto1} can be compared in Figure~\ref{fig:illpose}.

\begin{figure}[H]
	\begin{center}
		\includegraphics[width=6.2cm]{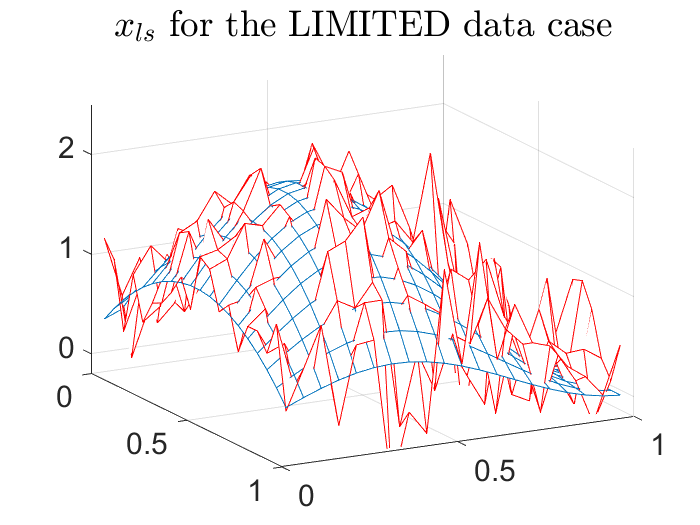}
		\includegraphics[width=6.2cm]{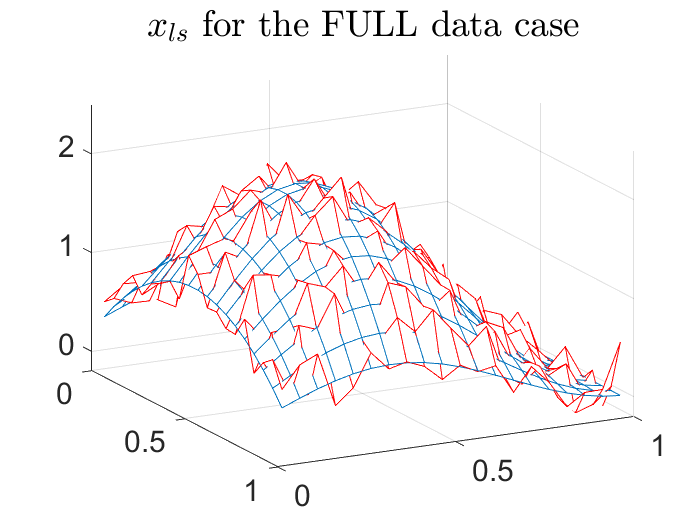}
		\caption{Phenomen of ill-posedness of deautoconvolution.}
		\label{fig:illpose}
	\end{center}
\end{figure}
As a consequence of the \textsl{ill-posedness phenomenon} of the deautoconvolution problem in 2D we see from the \Cref{fig:illpose} the occurrence of
\textsl{strong oscillations} in both cases. The oscillations, however, are much heavier for the limited data case (left-hand graph) than for the full data case (right-hand graph). The difference is particularly pronounced for
function values on the rear triangle half of the underlying unit square. The discretized $L^2$-norms of the deviation $x_{ls}-x^\dagger$ correspond to relative errors of $34.54\%$ (left) and $13.92\%$ (right).

\section{Regularization methods} \label{sec:regu}

\subsection{Tikhonov regularization and regularization parameter choices}\label{sec:tik}
As a first approach to overcome the ill-posedness of the two-dimensional deautoconvolution problem, we adopt the variational (Tikhonov-type) regularization, which is well-developed for solving ill-posed nonlinear operator equations. For this approach, the stable approximate solutions (\textit{regularized solutions}) $\xad$ are global minimizers of the optimization problem
\begin{equation}\label{eq:Tikreg}
	T_\alpha^\delta(x):=\|F(x)-y^\delta\|_Y^2+\alpha\, \mathcal{R}(x)\to \min_{x\in \mathcal{D}(F) \subseteq X}\,,
\end{equation}
with \textsl{regularization parameter} $\alpha>0$ and some  \textsl{penalty functional} $\mathcal{R}:X \to [0,+\infty]$ possessing the domain $\mathcal{D}(\mathcal{R}):=\{x \in X=L^2([0,1]^2): \mathcal{R}(x)< \infty\}$.
The penalty functional is assumed to be \textsl{stabilizing, convex} and \textsl{weakly sequentially continuous} such that for the autoconvolution operator $F$, which is \textsl{weakly sequentially closed}, the general theory of variational regularization (see, e.g, \cite[Section~4.1]{Schusterbuch12} and \cite{EHN96,Scherzetal09}) with respect to \textsl{existence, stability} and \textsl{convergence} of the Tikhonov-regularized solutions
$\xad \in \mathcal{D}(F) \cap \mathcal{D}(\mathcal{R})$ applies. The following three penalty functionals are under consideration in this study:
\begin{itemize}
	\item Classical norm square penalty
	\[\mathcal{R}_1(x):=\|x-\bar{x}\|_X^2\,,\]
	with prescribed reference element $\bar{x}\in X$ and $\mathcal{D}(\mathcal{R}_1)=X$. Notably, $\bar{x}\in X$ expresses some a-priori knowledge about the potential solution.
	\item Gradient norm square penalty \[\mathcal{R}_2(x):=\int_{[0,1]^2}\|\nabla x\|^2_2 \,dt_1 dt_2\] where $\nabla x=(\frac{\partial x}{\partial t_1},\frac{\partial x}{\partial t_2})$ denotes the gradient with respect to both variables $t_1, t_2$ and $\|\cdot\|_2$ is the Euclidean norm. Here, we have $\mathcal{D}(\mathcal{R}_2)=H^1([0,1]^2)$. For this setting, the solution is assumed to have a certain smoothness.
	\item Total variation penalty
	\[
	\mathcal{R}_3(x):= \|x\|_{TV([0,1]^2)} = \int_{[0,1]^2}\|\nabla x\|_1 \,dt_1 dt_2,
	\]
	where $\mathcal{D}(\mathcal{R}_3)=BV([0,1]^2):=\{x\in L^1([0,1]^2): \,\|x\|_{TV([0,1]^2)}<\infty\}$ is the space of bounded variation over the unit square $[0,1]^2$. This approach is originally introduced for image restoration with the special aim of retaining the information on edges in an image, i.e., the penalty should work for solutions possessing jumps. A detailed analysis of TV-regularization can, for example, be found in \cite{AcarVogel94,Burger13}.
\end{itemize}

\smallskip

In a first step, we want to choose the optimal regularization parameter $\alpha_{opt}$ for each input noise level $\delta$ and corresponding $y^\delta$ according to
\begin{equation}\label{eq:opt}
	\alpha_{opt}(\delta)=\text{argmin}_{\alpha>0}\|x_\alpha^\delta-x^\dagger\|_X.
\end{equation}
It is well-known that for the practical use of \textsl{a priori choices} for finding the regularization parameter $\alpha=\alpha(\delta)$, some smoothness information about the exact solution $x^\dagger$ is required,
which is normally not available. Therefore, \textsl{a posteriori choices} $\alpha=\alpha(\delta,y^\delta)$ are exploiting the measured noisy data $y^\delta$ in combination with knowledge of the noise level $\delta>0$  are
an appropriate alternative (cf.~\cite[Sect.~3.1]{EHN96}).  Under the limit conditions
\begin{equation}\label{eq:limitcond}
	\alpha(\delta,y^\delta)\to 0\quad\text{and}\quad \frac{\delta^2}{\alpha(\delta,y^\delta)}\to 0 \;\text{ as }\; \delta\to 0
\end{equation}
the regularized solutions $x_\alpha^\delta$ solving the optimization problem \eqref{eq:Tikreg} may possess a subsequence which converges to a exact solution $x^\dagger$ as $\delta\to 0$.
Due to the Fr\'echet differentiability of the autoconvolution operator $F$ and convexity of the applying penalty functionals $\mathcal{R}(x)$, for the numerical experiment below, we implement as second step the \emph{sequential discrepancy principle (SDP)} which was analyzed, for example, in \cite{Anzengruber14}.
\begin{definition}
	For given $\tau>1$, $\alpha_0>0$, $0<q<1$, a parameter $\alpha_{SDP}$ is chosen from the set $\Delta_q:=\{\alpha_l:\,\alpha_l=q^l\alpha_0,\,l\in \mathbb{Z}\}$ according to the sequential discrepancy principle (SDP), if
	\begin{equation}\label{eq:SDP}
		\|F(x_{\alpha_{SDP}}^\delta)-y^\delta\|_Y\leq \tau\delta<\|F(x_{\alpha_{SDP}/q}^\delta)-y^\delta\|_Y
	\end{equation}
	holds true.
\end{definition}

We can directly apply Theorem 1 from \cite{Anzengruber14} to our autoconvolution problem and conclude that, with some $\bar{\delta}>0$, the regularization parameters $\alpha_{SDP}=\alpha_{SDP}(\delta,y^\delta)$ chosen according to SDP exist for $0<\delta<\bar{\delta}$ and satisfy the limit conditions \eqref{eq:limitcond}. Then the associated regularized solutions $x_{\alpha_{SDP}(\delta,y^\delta)}^\delta$ converge (at least in the sense of subsequences) to exact solutions $x^\dagger$ as $\delta\to 0$  and moreover $\lim \limits_{\delta\to 0} \mathcal{R}(x_{\alpha_{SDP}(\delta,y^\delta)}^\delta)=\mathcal{R}(x^\dagger)$ holds true.

In a third step we search for \textsl{heuristic choices}  of the regularization parameter $\alpha>0$ provided that the noise level $\delta>0$ is not available or reliable. In our numerical case studies we focus only on the quasi-optimality criterion to find  $\alpha_{qo}=\alpha_{qo}(y^\delta)$, see \cite{Bakushinskii84,BauerKin08,KinNeu08} and references therein.

\begin{definition}
	For  sufficiently large $\alpha_0>0$ and for some $0<q<1$, we call the parameter $\alpha_{qo}$ chosen from the set $\Delta_q:=\{\alpha_l: \,\alpha_l=q^l\alpha_0,\,l\in \mathbb{N}\}$ according to
	\begin{equation}\label{eq:qo}
		\alpha_{qo}(y^\delta)=\argmin_{\alpha_l \in \Delta_q}\|x_{\alpha_l}^\delta-x_{\alpha_{l+1}}^\delta\|_X
	\end{equation}
	quasi-optimal regularization parameter.
\end{definition}

\subsection{An iteratively regularized Gauss-Newton method}

As an alternative we can also consider the iteratively regularized Gauss-Newton method (IRGNM) and find the minimizers of the functional
\begin{equation}\label{eq:irgnm}
	J_{\alpha_n}^\delta(x):=\|F(x_n^\delta)+F'(x_n^\delta)(x-x_n^\delta)-y^\delta\|_Y^2+{\alpha_n} \mathcal{R}(x) \to \min_{x\in \mathcal{D}(F) \subseteq X},
\end{equation}
with some initial guess $x_0^\delta \in X$ for a fixed noise level $\delta>0$. Therein,  $F'(x): X \to Y$ is the Fr\'echet derivative of $F$  at $x\in X$ and a sequence $\left(\alpha_n\right)_{n \in \N}$ of regularization parameters satisfies
\begin{equation}\label{eq:irgnm_regpar}
	1 \leq \frac{\alpha_{n}}{\alpha_{n+1}} \leq C
\end{equation}
with some constant $C > 0$. The central advantage of \eqref{eq:irgnm} over \eqref{eq:Tikreg} is that $x_{n+1}^\delta$ is defined as the solution of a (due to linearity of $F'(x_n^\delta)$) \textbf{convex} optimization problem, which can efficienlty be tackled by algorithms such as Chambolle-Pock \cite{ChamPock11} or FISTA \cite{BeckTeb14}.
By use of norm square penalty $\mathcal{R}_1(x)$,
\eqref{eq:irgnm} can be solved explicitly as
\begin{equation}\label{eq:irgR1}
	\left(F' (x_n^\delta)^*\left[F'(x_n^\delta)\right] + \alpha_n \mathbb{I} \right) x =  F' (x_n^\delta)^*\left[F'(x_n^\delta)x_n^\delta+y^\delta - F(x_n^\delta)\right] + \alpha_n \bar x.
\end{equation}

Since the gradient operator is also linear, we can solve the linear equation
\begin{equation}\label{eq:irgR2}
	\left(F'(x_n^\delta)^\star\left[ F'(x_n^\delta)\right]+\alpha_n\nabla^\star\nabla \right) x= F'(x_n^\delta)^\star\left[F'(x_n^\delta)x_n^\delta+y^\delta-F(x_n^\delta)\right]
\end{equation}
for the use of gradient norm square penalty $\mathcal{R}_2(x)$.

As a computational drawback, however, a full sequence of minimization problems (or linear equations) has to be solved. A convegence analysis for the IRGNM as depicted in \eqref{eq:irgnm} can e.g. be found in \cite{HohWer14, Wer15}.

Note that a similar approach has been proposed in \cite{Schleicher83}, where the least-squares residuum $\norm{F(x) - y^\delta}{Y}^2$ is minimized on a finite dimensional ansatz space $x \in \text{span} \left\{\mu_1, ..., \mu_n\right\}$ (e.g. consisting of splines) by linearization and iterative updating. In contrast to \eqref{eq:Tikreg} and \eqref{eq:irgnm}, regularization is there obtained by restriction to a finite-dimensional space, but the computational procedure is besides this comparable to our update formula in \eqref{eq:irgnm}.

Instead of choosing a regularization parameter $\alpha$ in Tikhonov regularization, here we have to select an appropriate stopping index $n \in\N_0$. This can in principle be done by the same rules as discussed in section \ref{sec:tik}.  The running index $n\in \mathbb{N}_0$ can be selected
\begin{itemize}
	\item in the best case as $n_{opt}$ with \[n_{opt}(\delta)= \argmin_{n \in \N_0}\norm{x_n^\delta -x^\dagger}{X};\]
	\item according to a-posteriori sequential discrepancy principle for a given constant parameter $\tau > 1$ as $n_{SDP}$ if \[
	\norm{F(x_{n_{SDP}}^\delta) - y^\delta}{Y} \leq \tau \delta<\norm{F(x_{n_{SDP}-1}^\delta) - y^\delta}{Y}
	\] holds true, see \cite{kns08} and the references therein.
\end{itemize}

\section{Numerical treatment} \label{sec:numerics}
\subsection{Discretization via the composite midpoint rule}

To discretize the continuous problem, we consider two different approaches. The first option is to divide each direction of the unit square equidistantly in $n$ partitions with the uniform length $h:=1/n$. In order to discretize the nonlinear convolution equation or deduce the discretized forward operator, it is reasonable to replace the function values $x(t_1,t_2)$ and $y(s_1,s_2)$ by countable values $x_{i,j}$ and $y_{k,l}$ with
\[x_{i,j}:=x(\tfrac12(i+(i-1))h,\tfrac12(j+(j-1))h),\qquad y_{k,l}=y(kh,lh)\] for all $i,j=1,\dots,n$ and $k,l=1,\dots,n$, respectively.

\subsubsection{Discretization of forward operator for the limited data case}

The autoconvolution equation in the limited data case can then by means of the composite midpoint rule be approximated by the discrete equation
\begin{equation}\label{eq:autonum01}
	\sum_{j=1}^l\sum_{i=1}^kh^2x_{k-i+1,l-j+1}x_{i,j}=y_{k,l}.
\end{equation}
Since the function $x$ takes value only for $0\leq t_1\leq s_1$ and $0\leq t_2\leq s_2$, only the indices $i\leq k$ and $j\leq l$ need to be taken into consideration in the discretized version accordingly.
Transforming the values $x_{i,j}$ and $y_{k,l}$ into vectors $\underline{x}:=(x_1,\dots,x_p,\dots,x_{n^2})^\top$ and $\underline{y}:=(y_1,\dots,y_q,\dots,y_{n^2})^\top$, respectively, with $p:=(i-1)\cdot n+j$ and $q:=(k-1)\cdot n+l$ for  all $p,q=1,\dots,n^2$, we can rewrite the weakly nonlinear forward operator of autoconvolution as
\begin{equation}\label{eq:ForwardF1}
	{F}_1(\underline{x}):=h^2M_1(\underline{x})\underline{x},
\end{equation}
where $M_1(\underline{x})\in \mathbb{R}^{n^2\times n^2}$ is a lower triangle block matrix and has the  structure
\begin{small}
	\[M_1(\underline{x}):=\begin{pmatrix}
		B_1&0&\cdots&0\\
		B_2&B_1&\cdots&0\\
		\vdots&\ddots&\ddots&\vdots\\
		B_n&B_{n-1}&\cdots&B_1
	\end{pmatrix} \text{ with }
	B_m=\begin{pmatrix}
		x_{(m-1)n+1}&0&\cdots&0\\
		x_{(m-1)n+2}&x_{(m-1)n+1}&\ddots&0\\
		\vdots&\ddots&\ddots&\vdots\\
		x_{mn}&x_{mn-1}&\cdots&x_{(m-1)n+1}
	\end{pmatrix}
	\]
\end{small}
for $1 \leq m \leq n$. The first derivative of ${F}_1(\underline{x})$ can then be easily obtained as
\[{F}_1'(\underline{x})=2h^2M_1(\underline{x}).\]

\subsubsection{Discretization of forward operator for the full data case}

In the full data case, the discrete forward operator of the autoconvolution equation can be derived in a similar way. On four subareas $\left[0,1\right]^2$, $\left[0,1\right]\times\left(1,2\right)$, $\left(1,2\right)\times\left[0,1\right]$ and $\left(1,2\right)^2$ we have the discrete equations as:
\begin{equation}\label{eq:autonum02}
	\begin{aligned}
		&\sum_{j=1}^l\sum_{i=1}^kh^2x_{k-i+1,l-j+1}x_{i,j}=y_{k,l}&&\quad \text{for}\quad 1\leq k,l\leq n\\
		&\sum_{j=l-n}^{n-1}\sum_{i=1}^kh^2x_{k-i+1,l-j}x_{i,j+1}=y_{k,l}&&\quad \text{for}\quad 1\leq k\leq n \text{  and  } n+1\leq l\leq 2n-1\\
		&\sum_{j=1}^{l}\sum_{i=k-n}^{n-1}h^2x_{k-i,l-j+1}x_{i+1,j}=y_{k,l}&&\quad \text{for}\quad n+1\leq k\leq 2n-1 \text{  and  } 1\leq l\leq n\\
		&\sum_{j=l-n}^{n-1}\sum_{i=k-n}^{n-1}h^2x_{k-i,l-j}x_{i+1,j+1}=y_{k,l}&&\quad \text{for}\quad n+1\leq k,l\leq 2n-1. \\
	\end{aligned}
\end{equation}
On the outside boundary of $[0,2]\times [0,2]$, i.e. for either $k=2n$ or $l=2n$, $y$ vanishes. Note that the grid width remains $h=1/n$. The discretized forward operator can be written as
\[F_2(\underline{x}):=h^2 M_2(\underline{x})\underline{x}\]
with an extended block matrix ${M}_2(\underline{x})\in \mathbb{R}^{4n^2\times n^2}$, where  \\
\begin{footnotesize}
	${M}_2(\underline{x}):=\begin{pmatrix}
		B_1&0&\cdots&\dots&0\\
		C_1&0&\cdots&\dots&0\\
		B_2&B_1&\cdots&\dots&0\\
		C_2&C_1&\cdots&\dots&0\\
		\vdots&\ddots&\ddots&\ddots&\vdots\\
		B_n&B_{n-1}&\cdots&\cdots&B_1\\
		C_n&C_{n-1}&\cdots&\cdots&C_1\\
		0&B_n&B_{n-1}&\cdots&B_{2}\\
		0&C_n&C_{n-1}&\cdots&C_{2}\\
		\vdots&\ddots&\ddots&\ddots&\vdots\\
		0&0&0&\cdots&B_n\\
		0&0&0&\cdots&C_n\\
		0&0&0&\cdots&0\\
		0&0&0&\cdots&0\\
	\end{pmatrix}$ \end{footnotesize}
with  \begin{footnotesize}
	$C_m=\begin{pmatrix}
		0&x_{mn}&x_{mn-1}&\cdots&x_{(m-1)n+2}\\
		0&0&x_{mn}&\ddots&x_{(m-1)n+3}\\
		\vdots&\ddots&\ddots&\vdots\\
		0&0&0&\cdots&x_{mn}\\
		0&0&0&\cdots&0
	\end{pmatrix}$ \end{footnotesize} and
$B_m$ as stated above for $1\leq m \leq n$.

The first derivative of the extended forward operator is then given as $ F_2'(\underline{x})=2h^2 M_2(\underline{x}).$

\subsubsection{Discretization of the penalty functionals}

The penalty functionals $\mathcal{R}_1$ and $\mathcal{R}_2$ can be discretized straight-forward, while the derivatives of $\nabla x$ are approximated by using finite differences. However, in order to ensure differentiability of the discretization of $\mathcal R_3$, we take the approximation $\sqrt{|w_1|^2 + |w_2|^2 + \beta}$ for the $1$-norm $\| w \|_1 = |w_1| + |w_2|$ of a vector $w = (w_1,w_2)$ with a small positive parameter $\beta\in (0,1)$. This idea leads to the discretization
\begin{equation}\label{eq:TVnum}
	\begin{aligned}
		\mathcal{R}_3(x)&\approx h\left(\sum_{i=1}^{n-1}\sum_{j=1}^{n-1}\sqrt{(x_{i+1,j}-x_{i,j})^2+(x_{i,j+1}-x_{i,j})^2+h^2\beta^2}\right.\\
		&\left.+\sum_{i=1}^{n-1}\sqrt{(x_{i+1,n}-x_{i,n})^2+h^2\beta^2}+ \sum_{j=1}^{n-1}\sqrt{(x_{n,j+1}-x_{n,j})^2+h^2\beta^2}\right),
	\end{aligned}
\end{equation}

\subsection{Discretization of forward operator via the Fourier transform}

If we assume that the function $x\in L^2((0,1)^2)$ can be extended on the whole domain $\mathbb{R}^2$ with the support of $x$ $\supp(x)\subset [0,1]\times [0,1]$, i.e. $x(t)=0$ for $t=(t_1,t_2)\notin [0,1]^2$, then we have the Fourier transform of $x$  as
\begin{equation}\label{eq:FT}
	\mathcal{F}(x)(\omega):=\frac{1}{\sqrt{2\pi}}\int_{[0,1]^2}x(t)e^{i\omega t}{d}t,\quad \omega\in \mathbb{R}^2.
\end{equation}
According to the convolution theorem, the autoconvolution operator can be represented by
\begin{equation}\label{eq:AutoFT}
	F(x)=\mathcal{F}^{-1}(\mathcal{F}(x)^2))
\end{equation}
Since the Fourier transform operator $\mathcal{F}$ w.r.t. $x$ is linear, we can easily obtain the Fr\'echet derivative and its adjoint operator of autoconvolution forward operator
\begin{equation}\label{eqs:AutoDerFT}
	\begin{aligned}
		&F'(x)(u)=\mathcal{F}^{-1}(2\mathcal{F}(x)\mathcal{F}(u))\\
		&F'(x)^\star(v)=\mathcal{F}^{-1}(2\overline{\mathcal{F}(x)}\mathcal{F}(v))\\
	\end{aligned}
\end{equation}
for $u\in X$, $v\in Y$ and $\overline{\cdot}$ denotes the conjugated complex value.

A discretization of the above formulas is directly available by means of the Fast Fourier Transform (FFT) and its inverse (IFFT). However, we should take into account that they consider \textit{periodic} functions, and hence the corresponding discretization will be erroneous especially close to the boundary and is not able to distinguish between the limited and the full data case.

For the limited data case we perform a zero-padding , i.e. replace the discretization $\underline{x} \in \R^{n \times n}$ by an extended matrix $\underline{x}^z \in \R^{2n \times 2n}$ of the form
\[
\underline{x}^z = \begin{pmatrix}
	\underline{x}  & 0 \\0  & 0
\end{pmatrix}.
\]
More precisely, if we denote the zero-padding operator by $Z : \R^{n\times n}\to \R^{2n\times 2n}$ and the corresponding left-inverse (restriction) by $R : \R^{2n \times 2n} \to \R^{n \times n}$, then we obtain the discretization of \eqref{eq:AutoFT} as
\begin{equation}\label{eqs:AutoFTdiscretization}	
	F(x) \approx R\left(\text{IFFT} \left(\text{FFT}(Z(\underline x))^2\right)\right).
\end{equation}
and the corresponding discretizations of \eqref{eqs:AutoDerFT}  as
\begin{equation}
	\begin{aligned}
		&F'(x)(u)\approx R\left(\text{IFFT}(2\text{FFT}(Z(\underline x))\text{FFT}(Z(\underline u)))\right)\\
		&F'(x)^\star(v)\approx R\left(\text{IFFT}(2\overline{\text{FFT}(Z(\underline x))}\text{FFT}(Z(\underline v)))\right).\\
	\end{aligned}
\end{equation}

The only difference between the limited data and the full data case is  whether the restriction $R$ as to be applied as a very last step or not.

\subsection{Computational implementation}

In order to tackle the autoconvolution problem in a stable way, we can either solve the Tikhonov regularized problem \eqref{eq:Tikreg} with diverse penalty functionals or solve the iteratively regularized problem \eqref{eq:irgnm} equipped with norm square penalty $\mathcal{R}_1(x)$ or gradient norm square penalty $\mathcal{R}_2(x)$. For both approaches, we need to initialize the regularization parameter $\alpha_0$ and set the iteration step $q\in (0,1)$ with $\alpha_l=q^{l}\alpha_0$ for $l=1,2,3,\dots$.

To solve the Tikhonov regularized problem \eqref{eq:Tikreg} with $\alpha:=\alpha_l$, we  use the discretization via the composite midpoint rule and consider its first-order optimality condition
\begin{equation}\label{eq:disc_Tik}
	2h^2(M(\underline{x}))^\top(h^2M(\underline{x})\underline{x}-y^\delta)+\alpha \mathcal{R}'(\underline{x})=0,
\end{equation}
which will be solved by using a damped Newton method. In this nonlinear equation, either $M:=M_1$ is in use for the limited data case or $M:=M_2$ for $y^\delta$ in the full data case. Note that the Newton-type method for solving \eqref{eq:disc_Tik} is also more or less a iterative procedure and  needs initialization of solution as well.

For our academical experiments ($x^\dagger$ is known), we carry out the process according to the following conceptional algorithm:
\begin{algorithm}[H]
	\caption{Conceptional algorithm for solving \eqref{eq:Tikreg}}
	\label{alg:Tikalg}
	\begin{algorithmic}
		\State \textbf{S0:} Let $\{\delta_k\}_{k\in \mathbb{N}}$ be a finite sequence of positive noise levels tending to $0$ as $k\to\infty$. Fix $\alpha_0>0$ and $0<q<1$. Set $x_0$ be a starting point, $l:=0$,  $k:=0$ and $x_{\alpha_l}^{\delta_k}:=x_0$. \\
		\textbf{S1:}  Compute the error $E_{l}^k=\|x_{\alpha_l}^{\delta_k}-x^\dagger\|_X$.\\
		\textbf{S2:} Solve the discretized problem \eqref{eq:disc_Tik} for fixed  $\delta:=\delta_{k}$ and $\alpha:=\alpha_l$ using a damped Newton method with starting point $x_{\alpha_l}^{\delta_k}$. Let $x_{\alpha_{l+1}}^{\delta_k}$ be the associated solution and compute the error $E_{l+1}^k=\|x_{\alpha_{l+1}}^{\delta_k}-x^\dagger\|_X$.\\
		\textbf{S3:} If $E_{l+1}^k<E_{l}^k$, then $\alpha_{l+1}:=q\alpha_l$, $l:=l+1$ and go to \textbf{S2.} Otherwise save $x_{\alpha_l}^{\delta_k}$ and go to \textbf{S4}.\\
		\textbf{S4:} Set $x_{\alpha_0}^{\delta_{k+1}}:=x_{\alpha_l}^{\delta_k}$, $l:=0, k:=k+1$ and go to \textbf{S2}.
	\end{algorithmic}
\end{algorithm}

Note that this algorithm  can be efficiently implemented due to the square rate of convergence of Newton-type method.  However, the derivation  of the first and second derivative of the autoconvolution operator and penalty functionals costs us a certain effort.

To solve the iteratively regularized problem  \eqref{eq:irgnm} with either $\mathcal{R}:=\mathcal{R}_1$ or $\mathcal{R}:=\mathcal{R}_2$, more precisely, to solve the linear equations   \eqref{eq:irgR1} and \eqref{eq:irgR2}, we make use of Fourier transform technique to discretize the forward operator and apply the CG (conjugate gradient) method. Besides benefiting the fast computation time, we can also save the trouble of derivation of second derivative matrices of all functionals. Actually, the associated algorithm is similar to \Cref{alg:Tikalg}. The notation $x_{\alpha_l}^{\delta_k}$ in \Cref{alg:Tikalg} corresponds to the solution of  \eqref{eq:irgR1} and \eqref{eq:irgR2} on the step $n:=l$ for $\delta:=\delta_k$.

\section{Numerical experiments} \label{sec:studies}
In this section, we present case studies of numerical experiments based on both examples that had been introduced in Subsection~\ref{sec:example}. For all three penalties $\mathcal{R}_i(x)\;(i=1,2,3)$ introduced above, properties and
accuracy of regularized solutions are compared and illustrated for various  regularization parameter choice rules. We define a sequence of finite discrete noise levels $\{\delta_k\}_{k\in \mathbb{N}}$ with corresponding relative noise levels between $10\%$ and $0.05\%$. Let input data $y^{\delta_k}\in Y$ satisfy $\|y^{\delta_k}-y\|_Y=\delta_k$. In order to find appropriate regularization parameters $\alpha$ for each fixed $\delta_k$, we set for the  sequence $(\alpha_l)_{l\in \mathbb{N}}$
the starting value $\alpha_0=1$, step size $q=0.5$ and thus $\alpha_{l+1}=\alpha_l/2$, see \Cref{alg:Tikalg}. On the discretized domain with discretization level $n=20$, we solve the nonlinear equation  \eqref{eq:disc_Tik}, discretized counterparts of linear equations \eqref{eq:irgR1} and \eqref{eq:irgR2}  for each $\delta_k$ and the regularization parameter $\alpha_l$. Additionally, the constant function $x_0\equiv 1$ is set as initialization for all computations. For the use of penalty term $\mathcal{R}_1(x)$, the reference element is set as a constant function $\bar{x}=0.5$. If we set $\bar{x}=-0.5$, the solutions $x_{\alpha_l}^{\delta_k}$ will converge to $-x^\dagger$ for both examples. Following the SDP choosing rule to determine the regularization parameters $\alpha_{SDP}^{\delta}$ in Tikhonov regularized problem or the running index $n_{SDP}$ in iteratively regularized problem, the constant parameter $\tau=1.2$ is fixed for our experiments. By use of smoothed total variation the smoothing parameter $\beta$ is fixed as $0.1$ in the discretized TV penalty \eqref{eq:TVnum}.

\subsection{Results for Example 1}
\begin{figure}[H]
	\begin{center}
		\includegraphics[width=12cm]{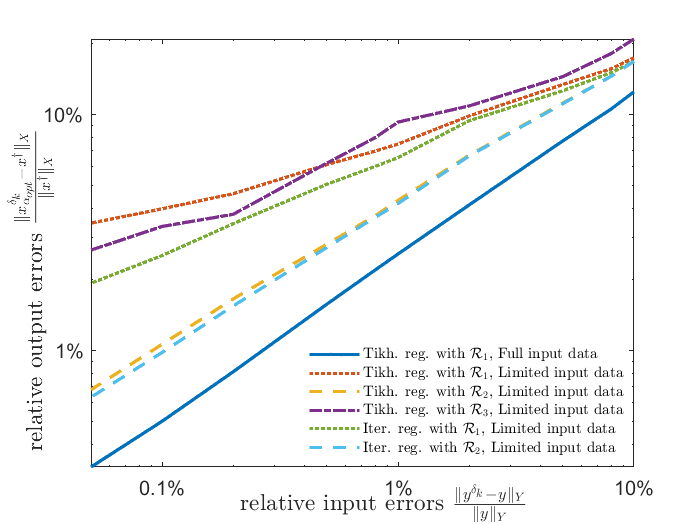}
		\caption{Comparison of relative  error norms of regularized solutions $x_{\alpha_{opt}}^{\delta_k}$ with optimal regularization parameter.}
		\label{fig:ex1_optsol}
	\end{center}
\end{figure}

First we compare in  \Cref{fig:ex1_optsol} the relative output errors of the regularized solutions $x_{\alpha_{opt}}^{\delta_k}$ obtained by using diverse penalties and by distinguishing the full data case and the limited data case. Since, for all penalties, the accuracies obtained with the full input data are uniformly better
than those obtained with limited input data, we have illustrated in the following figure for the full data case only accuracies to the norm square penalty $\mathcal{R}_1(x)$ by use of Tikhonov regularization. On the other hand, results for all penalties and both regularization methods are displayed for the limited data case.

Obviously, the operation with the full input data brings the smallest and best output error. For the limited data case, the gradient norm square penalty is  the most suitable for this example. A perspicuous reason is the high level of smoothness of the exact solution $x^\dagger$. The quality of Tikhonov-regularized solutions based on the classical norm square penalty $\mathcal{R}_1$ and of the TV-penalty $\mathcal{R}_3$ is almost indiscernible.
Moreover, when using this gradient penalty $\mathcal{R}_2$, the quality of results for Tikhonov regularization and iterative regularization are nearly not indistinguishable.

\smallskip

In \Cref{tab:ex1_conv} we present H\"older exponents $\kappa \in (0,1)$ estimated by regression from a series of $\delta$-values, which emulate numerically convergence rates results for regularized solutions with best possible regularization parameter $\alpha_{opt}$ for each $\delta$. Precisely, we have listed
the exponents $\kappa$ such that approximately $\|x_{\alpha_{opt}}^\delta-x^\dagger\|_X\sim \delta^\kappa$ as $\delta \to 0$.

\begin{table}[H]
	\caption{Estimated H\"older exponents $\kappa \in (0,1)$ for H\"older convergence rates $\|x_{\alpha_{opt}}^\delta-x^\dagger\|_X\sim \delta^\kappa$ as $\delta \to 0$}
	\label{tab:ex1_conv}
	\begin{center}
		\begin{tabular}{| c| c| c |c|c|c|c|}
			\hline
			\multirow{2}{*}{\parbox{3cm}{data situation \\penalty variation }} & \multicolumn{6}{|c|}{full data case \hspace{1.5cm} limited data case}\\
			&\multicolumn{3}{|c|}{$Y=L^2([0,2]^2)$}  &\multicolumn{3}{|c|}{$Y=L^2([0,1]^2)$}\\
			\hline
			H\"older exponent $\kappa$ &$\mathcal{R}_1(x)$&$\mathcal{R}_2(x)$&$\mathcal{R}_3(x)$&$\mathcal{R}_1(x)$&$\mathcal{R}_2(x)$&$\mathcal{R}_3(x)$\\
			\hline
			Tikhonov reg. & $0.6946$&$0.6638$&$0.4685$&$0.3118$&$0.6015$&$0.3919$\\
			\hline
			Iterative reg. & $0.7184$&$0.6936$&$-$&$0.4088$&$0.6183$&$-$\\
			\hline
		\end{tabular}
	\end{center}
\end{table}

For the limited data case, the gradient norm square penalty also delivers the largest H\"older rate exponent for the regularized solutions with best possible regularization parameter, which is consistent with the insights from  \Cref{fig:ex1_optsol} above. Partially, the rate exponents for the iterative regularization method seems to be higher than in case of Tikhonov regularization.

Next, we compare in \Cref{tab:ex1_output} the regularized solutions with  respect to different regularization parameter chosen rules, representatively for a fixed noise level $\delta$ of a relative input error $\frac{\|y^{\delta}-y\|_Y}{\|y\|_Y}=1\%$:

\begin{table}[H]
	\caption{Comparison of relative error norms of regularized solutions for a fixed relative input error $1\%$}
	\label{tab:ex1_output}
	\begin{center}
		\begin{tabular}{| c| c| c |c|c|c|c|}
			\hline
			relative input errors & \multicolumn{4}{|c|}{Tikhonov reg.}&\multicolumn{2}{|c|}{Iterative reg.}\\
			\hline
			&$Y=L^2([0,2]^2)$ &\multicolumn{5}{|c|}{ $Y=L^2([0,1]^2)$} \\
			&$\mathcal{R}_1(x)$&$\mathcal{R}_1(x)$&$\mathcal{R}_2(x)$&$\mathcal{R}_3(x)$&$\mathcal{R}_1(x)$&$\mathcal{R}_2(x)$\\
			\hline
			$\frac{\|x_{\alpha_{opt}}^{\delta}-x^\dagger\|_X}{\|x^\dagger\|_X}$&$2.56\%$&$7.49\%$&$4.32\%$&$9.28\%$&$6.56\%$&$4.19\%$\\
			\hline
			$\frac{\|x_{\alpha_{SDP}}^{\delta}-x^\dagger\|_X}{\|x^\dagger\|_X}$&$2.71\%$&$9.12\%$& $9.91\%$& $13.77\%$& $9.07\%$&$9.88\%$\\
			\hline
			$\frac{\|x_{\alpha_{qo}}^{\delta}-x^\dagger\|_X}{\|x^\dagger\|_X}$&$2.71\%$&$9.12\%$& $4.39\%$& $10.30\%$&- &-\\ 
			\hline
		\end{tabular}
	\end{center}
\end{table}
The regularized solutions with the regularization parameters via discrepancy principle can be obtained in a more stable way than the use of  quasi-optimality criterion, especially in pursuance of iteratively regularized Gauss-Newton-Method.

\begin{figure}[H]
	\begin{center}
		\includegraphics[width=6.2cm]{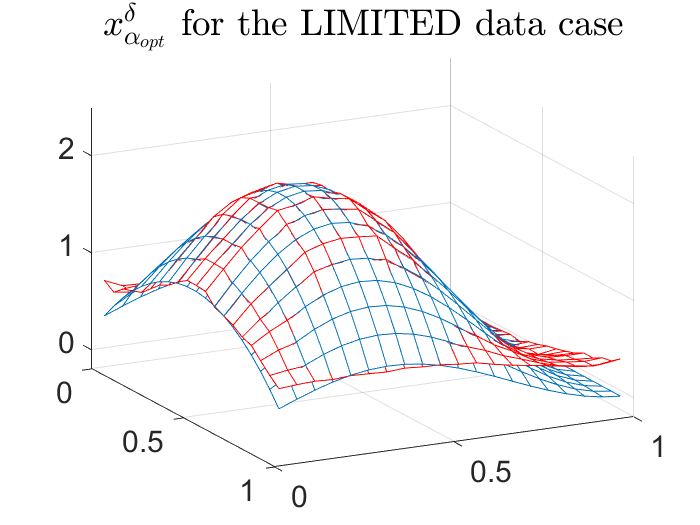}
		\includegraphics[width=6.2cm]{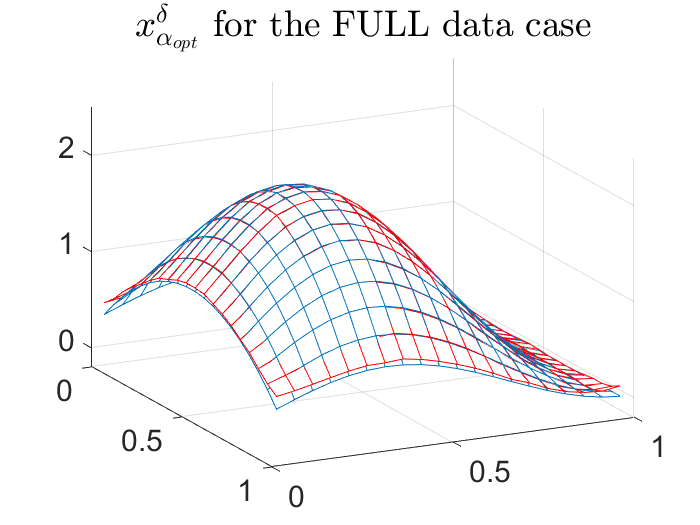}
		\caption{Regularized solutions with optimal regularization parameter for the limited and the full data case.}
		\label{fig:ex1_regsol}
	\end{center}
\end{figure}

For a more intuitive comparison of regularized solutions with least-square solutions which had been presented in \Cref{fig:illpose}, we now illustrate Tikhonov regularized solutions with best possible regularization parameter, obtained with classical norm square penalty $\mathcal{R}_1(x)$  for the same noise level of $0.8\%$ in \Cref{fig:ex1_regsol} in the limited and in the full data case, respectively.

There are almost no longer oscillations in the images above. However, we can still observe the deviation of $x_{\alpha_{opt}}^\delta$ to the exact $x^\dagger$ on the far back corner of the underlying unit square for the limited data case. These deviations can be eliminated by use of full input data.

\subsection{Results for Example 2}

Now we present analog results for Example~2 with $\xdag$ from  \eqref{eq:ex2}, which represents a non-smooth and non-factored function.

\begin{figure}[H]
	\begin{center}
		\includegraphics[width=12cm]{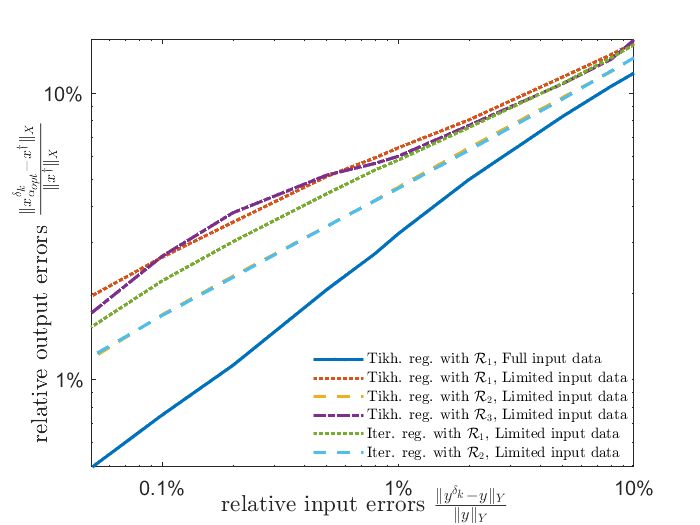}
		\caption{Comparison of relative  error norms of regularized solutions $x_{\alpha_{opt}}^{\delta_k}$ with optimal regularization parameter.}
		\label{fig:ex2_optsol}
	\end{center}
\end{figure}

The total variation penalty $\mathcal{R}_3(x)$  achieves progress among all penalties, since the function $x^\dagger$ possesses jumps and is much less smooth compared to Example~1.

\begin{figure}[H]
	\begin{center}
		\includegraphics[width=7.5cm]{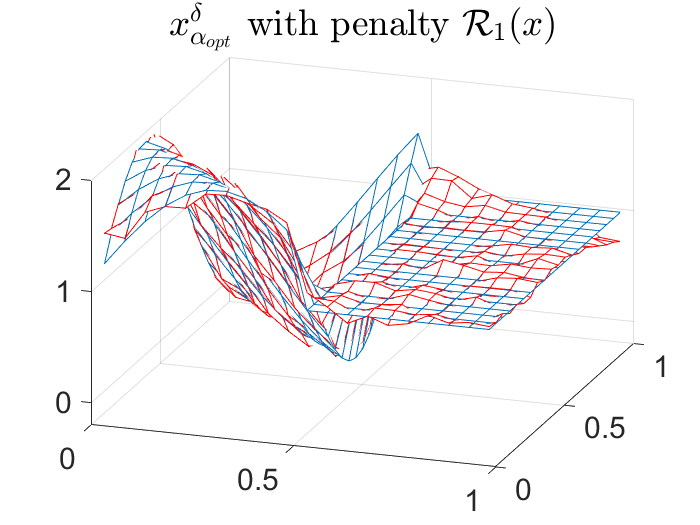}
		\includegraphics[width=7.5cm]{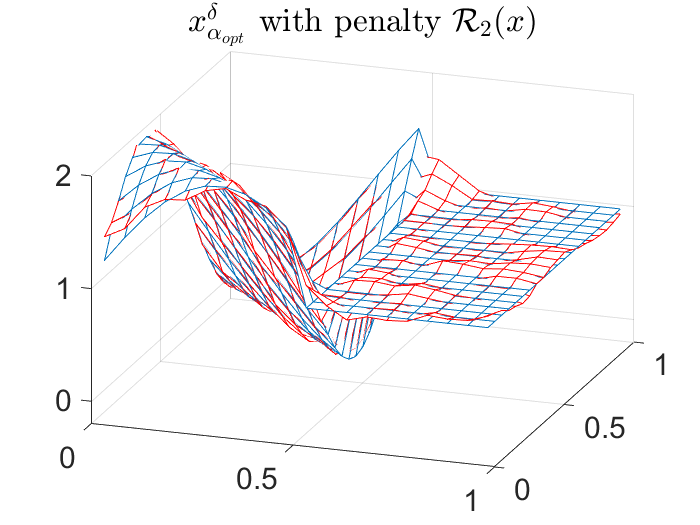}
		\includegraphics[width=7.5cm]{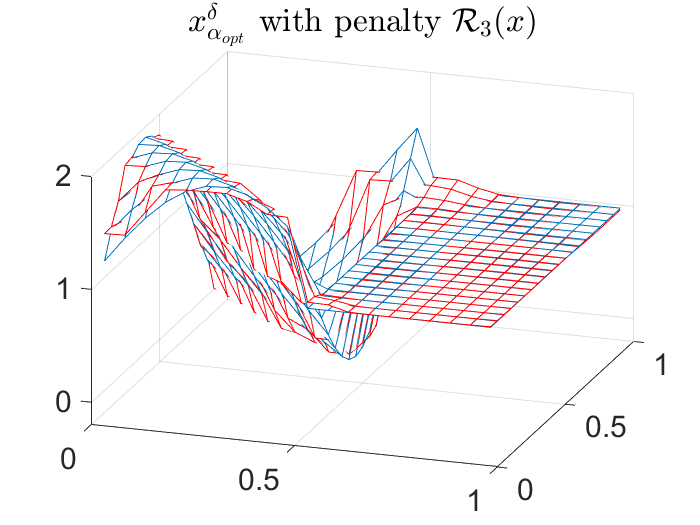}
		\caption{Regularized solutions obtained with optimal regularization parameter and for different penalties.}
		\label{fig:ex2_regsol}
	\end{center}
\end{figure}

The \Cref{fig:ex2_regsol} shows the shapes of Tikhonov-regularized solutions with best possible regularization parameters, calculated for different penalties and limited input data by a relative noise level of $0.8\%$. The involvement of TV penalty $\mathcal{R}_3(x)$  can be especially  observed on the rear area, where the function values tend to be constant. By contrast, the gradient norm square penalty $\mathcal{R}_2(x)$ makes the associated regularized solution more smooth on the front area.

We  again present also the estimated H\"older exponents $\kappa$ for regularized solutions $x_{\alpha_{opt}}^{\delta_k}$ with optimal regularization parameter.

\begin{table}[H]
	\caption{Estimated H\"older exponents $\kappa \in (0,1)$ for H\"older convergence rates $\|x_{\alpha_{opt}}^\delta-x^\dagger\|_X\sim \delta^\kappa$ as $\delta \to 0$}
	\label{tab:ex2_conv}
	\begin{center}
		\begin{tabular}{| c| c| c |c|c|c|c|}
			\hline
			\multirow{2}{*}{\parbox{3cm}{data situation \\penalty variation }} & \multicolumn{6}{|c|}{full data case \hspace{1.5cm} limited data case}\\
			&\multicolumn{3}{|c|}{$Y=L^2([0,2]^2)$}  &\multicolumn{3}{|c|}{$Y=L^2([0,1]^2)$}\\
			\hline
			H\"older exponent $\kappa$&$\mathcal{R}_1(x)$&$\mathcal{R}_2(x)$&$\mathcal{R}_3(x)$&$\mathcal{R}_1(x)$&$\mathcal{R}_2(x)$&$\mathcal{R}_3(x)$\\
			\hline
			Tikhonov reg. & $0.6059$&$0.6320$&$0.5083$&$0.3753$&$0.4522$&$0.3787$\\
			\hline
			Iterative reg. &$0.6699$&$0.6705$&-&$0.4164$&$0.4505$&-\\
			\hline
		\end{tabular}
		
	\end{center}
\end{table}

Regarding the  convergence rates of $x_{\alpha_{opt}}^\delta$ converging to $x^\dagger$ as $\delta\to 0$, for both penalties $\mathcal{R}_1$ and $\mathcal{R}_2$ the iterative regularization approach proves to be advantageous.

Similar to the Example 1, we list in \Cref{tab:ex2_output} the relative output errors of regularized solutions $x_{\alpha_{opt}}^{\delta}$, $x_{\alpha_{SDP}}^{\delta}$ and $x_{\alpha_{qo}}^{\delta}$ for the fixed noise level $\delta$ of relative input error $1\%$.

\begin{table}[H]
	\caption{Comparison of relative error norms of regularized solutions for a fixed relative input error $1\%$}
	\label{tab:ex2_output}
	\begin{center}
		\begin{tabular}{| c| c| c |c|c|c|c|}
			\hline
			relative input errors & \multicolumn{4}{|c|}{Tikhonov reg.}&\multicolumn{2}{|c|}{Iterative reg.}\\
			\hline
			&$Y=L^2([0,2]^2)$ &\multicolumn{5}{|c|}{ $Y=L^2([0,1]^2)$} \\
			&$\mathcal{R}_1(x)$&$\mathcal{R}_1(x)$&$\mathcal{R}_2(x)$&$\mathcal{R}_3(x)$&$\mathcal{R}_1(x)$&$\mathcal{R}_2(x)$\\
			\hline
			$\frac{\|x_{\alpha_{opt}}^{\delta}-x^\dagger\|_X}{\|x^\dagger\|_X}$&$3.22\%$&$6.43\%$&$4.68\%$&$6.01\%$&$5.83\%$&$4.64\%$\\
			\hline
			$\frac{\|x_{\alpha_{SDP}}^{\delta}-x^\dagger\|_X}{\|x^\dagger\|_X}$ & $5.75\%$&$8.59\%$& $ 7.75\%$&
			$10.60\%$& $9.06\%$& $8.92\%$\\
			\hline
			$\frac{\|x_{\alpha_{qo}}^{\delta}-x^\dagger\|_X}{\|x^\dagger\|_X}$&$3.47\%$&$7.23\%$& $ 5.83\%$&
			$6.01\%$& - &-\\
			\hline
		\end{tabular}
	\end{center}
\end{table}

By summarizing the case studies, we can state that it is always possible to solve the two-dimensional deautoconvolution problem in a rather stable way by either Tikhonov or iterative regularization. In this context, we
also obtain a reasonably accuracy when the choice of the regularization parameter is appropriate.

\section*{Acknowledgment} Yu Deng and Bernd Hofmann are supported by the German Science Foundation (DFG) under the grant~HO~1454/13-1 (Project No.~453804957).

\end{document}